\documentclass[12pt]{article}
\usepackage[table]{xcolor}
\usepackage{graphicx}
\usepackage{amsmath,amsfonts,amssymb}
\usepackage{graphicx}
\usepackage{multirow}
\usepackage{tikz,pgf}
\usepackage{url}
\usepackage{amsmath}
\usepackage{graphicx}
\usepackage{epsfig}
\usepackage{epstopdf}
\usepackage{color}
\usepackage{enumitem}
\def\disp{\displaystyle}

\def\tto{\;{\lower 1pt \hbox{$\rightarrow$}}\kern -10pt
\hbox{\raise 2pt \hbox{$\rightarrow$}}\;}
\def\Hat{\widehat}
\def\hat{\widehat}
\def\vt{\vartheta}

\def\Bar{\overline}
\def\ra{\rangle}
\def\la{\langle}
\def\kk{\kappa}
\def\ve{\varepsilon}
\def\epsilon{\varepsilon}
\def\B{\Bbb B}
\def\h{\hfill\Box}
\def\R{\Bbb R}

\def\ox{\bar{x}}

\def\oy{\bar{y}}
\def\oz{\bar{z}}

\def\op{\bar{p}}

\def\gph{\mbox{\rm gph}\,}

\def\cl{\mbox{\rm cl}\,}

\def\h{\hfill\square}

\def\ph{\varphi}
\def\emp{\emptyset}
\def\st{\stackrel}
\def\oR{\Bar{\R}}
\def\lm{\lambda}

\def\dd{\delta}
\def\al{\alpha}
\def\bb{\beta}
\def\Th{\Theta}

\def\ph{\varphi}
\def\emp{\emptyset}
\def\st{\stackrel}
\def\oR{\Bar{\R}}
\def\lm{\lambda}

\def\dd{\delta}
\def\al{\alpha}
\def\Th{\Theta}

\setlist[enumerate,1]{itemsep=0.0ex,parsep=0.5ex,label={\rm(\alph*)},leftmargin=*, align=left}
\newcounter{lk}

\usepackage[colorlinks=true]{hyperref}

\begin{document}
\begin{center}
{\bf\large Coincidence Points of Parameterized Generalized Equations with Applications to Optimal Value Functions}\\[2ex]
{\sc Aram V. Arutyunov}\footnote{V. A. Trapeznikov Institute of Control Sciences of Russian Academy of Sciences, Russia (arutyunov@cs.msu.ru). Research of this author was partly supported by the Russian Science Foundation, project 22-21-00863.}, {\sc Boris S. Mordukhovich}\footnote{Department of Mathematics, Wayne State University, USA (aa1086@wayne.edu). Research of this author was partly supported by the USA National Science Foundation under grant DMS-1808978, by the Australian Research Council under grant DP-190100555, and by Project 111 of China under grant D21024.} and {\sc Sergey E. Zhukovskiy}\footnote{V. A. Trapeznikov Institute of Control Sciences of Russian Academy of Sciences, Russia (s-e-zhuk@yandex.ru). Research of this author was partly supported by the Russian Science Foundation, project 20-11-20131.} 
\end{center}
\small{\bf Abstract.} The paper studies coincidence points of parameterized set-valued mappings (multifunctions), which provide an extended framework to cover several important topics in variational analysis and optimization that include the existence of solutions of parameterized generalized equations, implicit function and fixed-point theorems, optimal value functions in parametric optimization, etc. Using the advanced machinery of variational analysis and generalized differentiation that furnishes complete characterizations of well-posedness properties of multifunctions, we establish a general theorem ensuring the existence of parameter-dependent coincidence point mappings with explicit error bounds for parameterized multifunctions between infinite-dimensional spaces. The obtained major result yields a new implicit function theorem and allows us to derive efficient conditions for semicontinuity and continuity of optimal value functions associated with parametric minimization problems subject to constraints governed by parameterized generalized equations. 
\\[1ex]

\newtheorem{Theorem}{Theorem}[section]
\newtheorem{Proposition}[Theorem]{Proposition}
\newtheorem{Remark}[Theorem]{Remark}
\newtheorem{Lemma}[Theorem]{Lemma}
\newtheorem{Corollary}[Theorem]{Corollary}
\newtheorem{Definition}[Theorem]{Definition}
\newtheorem{Example}[Theorem]{Example}
\renewcommand{\theequation}{\thesection.\arabic{equation}}
\normalsize

\section{Introduction and Initial Discussions}\label{intro}
\setcounter{equation}{0}

The main attention of this paper is paid to studying the class of {\em parameterized generalized equations} given by 
\begin{equation}\label{pge}
0\in G(x,p)+F(x),
\end{equation}
where both mappings $G\colon X\times P\tto Y$ and $F\colon X\tto Y$ may be set-valued, i.e., being multifunctions. Although all the obtained results are new even in finite dimensions, we proceed in what follows in a broad generality; see the precise assumptions below. The generalized equation formalism was introduced by Robinson \cite{rob} in the case when the parameter-dependent mapping $G(x,p)$ in \eqref{pge} is single-valued, while the multifunction $F(x)$ is the normal cone mapping to a convex set. Robinson used his model to describe variational inequalities and  complementarity systems in optimization in the form of ``generalized equations" in order to obtain for them theoretical and numerical results similar to those for usual equations. Over the years, the original generalized equation formalism has been extended to different settings (in particular, when $F$ is a subdifferential of an extended-real-valued function), which cover a variety of models that naturally arise in optimization theory, numerical methods, and applications; see, e.g., the books \cite{bs,dr,is,m06} with the vast bibliographies therein. 

Nevertheless, the standard framework of \eqref{pge} with a single-valued mapping $G\colon X\times P\to Y$ does not encompass some classes of variational systems important for optimization theory and applications, which thus require the extended version of \eqref{pge}, where both mappings $G$ and $F$ are set-valued. Among such systems, we particularly mention {\em set-valued variational inequalities} initiated in \cite{ah} and {\em multiobjective optimization problems with equilibrium constraints} formulated and studied in \cite{bgm}. Furthermore, the parameterized generalized equations of type (\ref{pge}), described by both set-valued mappings $G$ and $F$, unavoidably arise in variational systems generated by first-order necessary optimality conditions for problems of {\em parametric nonsmooth optimization}:
\begin{equation}\label{nop}
\mbox{minimize }\;\ph(x,p)\;\mbox{ subject to }\;x\in\Omega,
\end{equation} 
where $\ph\colon X\times P\to\oR:=(-\infty,\infty]$ is an extended-real-valued lower semicontinuous function, and where $\Omega$ is a closed set. It is well known that necessary optimality conditions for \eqref{nop} are given in the form
\begin{equation}\label{noc}
0\in\partial_x\ph(x,p)+N(x;\Omega)
\end{equation}
in terms of certain subdifferential and normal cone constructions in appropriate Banach spaces; see \cite{m06} for more details. Hence we are again in the setting of \eqref{pge} with both set-valued mappings $G$ and $F$. Note that problems of type \eqref{nop} can be treated as lower-level problems in {\em nonsmooth bilevel programming} \cite{d,m18}. If both $\ph(\cdot,p)$ and $\Omega$ are convex, then \eqref{nop} is {\em equivalent} to the parameterised generalized equation \eqref{pge}. 

The underlying issue in the study of parameterized generalized equations is the {\em existence} of their solutions. Known results in this direction are obtained in the standard setting when $G$ in \eqref{pge} is {\em single-valued}, and this assumption is essential in the proofs of solvability. Indeed, the essence of the solvability proofs for \eqref{pge} consists of the following; see, e.g., \cite[Section~2.2]{dr}. If a set-valued mapping $F$ is metrically regular, while a single-valued mapping $G$ is Lipschitz continuous in $x$ with a sufficiently small Lipschitz constant, then the sum $G+F$ is metrically regular in $x$, which ensures the existence of solutions to \eqref{pge} in this case. The imposed single-valuedness assumption is {\em essential} for such an approach to the solvability of \eqref{pge}, since without it the sum of metrically regular and Lipschitzian mappings may not be metrically regular.  

In this paper, we develop a different approach to the solvability of parameterized generalized equations, which allows us to establish the existence of a solution to  \eqref{pge} in the general case where both $G$ and $F$ are set-valued. Given multifunctions $F$ and $G$ as in \eqref{pge}, consider a parametric family of their {\em coincidence points} 
\begin{equation}\label{cp}
F(x)\cap G(x,p)\ne\emp
\end{equation}
and aim to find efficient conditions that ensure the existence of a single-valued mapping $x(p)$ defined on some open subset of $P$ such that
\begin{equation*}
F\big(x(p)\big)\cap G\big(x(p),p\big)\ne\emp.
\end{equation*}
It is easy to see that the extended parameterized generalized equation \eqref{pge} with both set-valued mappings $F$ and $G$ is equivalent to the 
the parameterized coincidence point problems of type \eqref{cp} with the replacement of $G$ by $-G$ therein.

A particular case of \eqref{cp} is the classical problem on {\em implicit functions} formulated as follows. Given $f\colon X\times P\to Y$, find a single-valued mapping $x(p)$ such that $f(x(p),p)\equiv 0$ on some open subset of $P$. The main result of this paper on coincidence points of \eqref{cp} leads us to a new implicit function theorem.

The rest of the paper is organized in the following way. Section~\ref{prel} presents some preliminaries from variational analysis and generalized differentiation broadly used in formulations and proofs of the new results. In Section~\ref{sec:coin} we establish a major result of the paper providing a general parameterized coincidence point theorem. Section~\ref{sec:pge} contains applications of this theorem to deriving new results on the existence of solutions to parameterized generalized equations and implicit functions with error bound estimates. In Section~\ref{sec:marg} we develop applications of the obtained results to the study of semicontinuity, continuity, and Lipschitz continuity properties of optimal value functions in parametric minimization problems with constraints governed by parameterized generalized equations. The final Section~\ref{conc} summarizes the major results of the paper and discusses several topics of our further research. 

\section{An Overview from Variational Analysis}\label{prel}
\setcounter{equation}{0}

In this section, we overview some needed fundamental notions of variational analysis along with their characterizations via generalized differentiation. The reader can find more details in \cite{m06} in infinite-dimensional spaces and in \cite{m18,rw} for specifications in finite dimensions. Throughout the paper, the used notation and terminology are standard in variational analysis and generalized differentiation; see, e.g., \cite{m06,rw}. Recall that, in the setting of normed spaces $Z$, the symbol $B_Z(z,r)$ stands for the closed ball centered at $z\in Z$ with radius $r>0$, while we use $\B_Z$ to denote the unit closed ball centered at the origin of the space in question if no confusion arises. 

We start with the fundamental property of set-valued mappings that is mainly exploited in the paper. Given a multifunction $\Phi\colon X\tto Y$ between normed spaces and given nonempty sets $U\subset X$ and $V\subset Y$, we say that $\Phi$ enjoys 
the {\em covering property} with modulus $\al>0$ (or it has the $\al$-covering property) on $U$ relative to $V$ if
\begin{equation}\label{cov}
\Phi(x)\cap V+\al r\B_Y\subset\Phi(x+r\B_X)\;\mbox{ whenever }\;x+r\B_X\subset U\;\mbox{ as }\;r>0.
\end{equation}
Given $\ox\in X$ and $\oy\in\Phi(\oy)$ and taking in \eqref{cov} the sets $U$ and $V$ as neighborhoods of 
$\ox$ and $\oy$, respectively, we speak about the {\em $\al$-covering} of $\Phi$ {\em around} $(\ox,\oy)$. Note that the latter local covering property of $\Phi$ is also known as linear openness, or openness at a linear rate, of $\Phi$ around 
$(\ox,\oy)$. 

It has been well recognized in variational analysis (see, e.g., \cite[Theorem~1.52]{m06}) that the covering property of $\Phi$ around $(\ox,\oy)$ is equivalent to the {\em metric regularity} of $\Phi$ around this point with the reciprocal modulus, meaning that
\begin{equation}\label{mr}
{\rm dist}\big(x;\Phi^{-1}(y)\big)\le\al^{-1}{\rm dist}\big(y;\Phi(x)\big)\;\mbox{ for all }\;x\in U,\;y\in V,
\end{equation}
where ${\rm dist}(\cdot;\Omega)$ stands for the distance function associated with the set $\Omega$, and where $\Phi^{-1}\colon Y\tto X$ is the inverse mapping of $\Phi$. Furthermore, the covering/metric regularity property of $\Phi$ happens to be equivalent to a certain Lipschitzian behavior of the inverse mapping $\Phi^{-1}$; see \cite[Theorem~1.49]{m06} for the precise statements. 

The Lipschitzian property corresponding to \eqref{cov} and \eqref{mr} for the inverse mappings is formulated as follows. Given nonempty subsets $U\subset X$ and $V\subset Y$, the multifunction $\Phi\colon X\tto Y$ between normed spaces is called {\em Lipschitz-like} on $U$ relative to $V$ with some modulus $\ell\ge 0$ if we have the inclusion
\begin{equation}\label{lip}
\Phi(x)\cap V\subset\Phi(u)+\ell\|x-u\|\B_Y\;\mbox{ for all }\;x,u\in U.
\end{equation}
The local counterpart of \eqref{lip}, which is also known as the pseudo-Lipschitz or Aubin property of $\Phi$ around $(\ox,\oy)$, is formulated accordingly. 

A great advantage of generalized differentiation in variational analysis is that appropriate generalized differential constructions for multifunctions allow us to {\em completely characterize} the aforementioned fundamental properties of set-valued mappings around the references points. There are two types of characterizations of these properties labeled as {\em pointbased} and {\em neighborhood} ones; see \cite[Chapter~4]{m06}. The pointbased characterizations are expressed in terms of {\em robust} (stable with respect to perturbations of the initial data) coderivatives defined at the point in question, while the neighborhood ones are obtained in terms of certain ``precoderivatives" involving  points near the given one. The obtained pointbased characterizations are preferable since they use only the reference point and the employed coderivatives, besides being robust, admit comprehensive {\em calculus rules} based on {\em variational/extremal principles} of variational analysis as developed in \cite{m06}. In finite dimensional spaces, the pointbased and neighborhood characterizations are equivalent, and thus we can use the pointbased ones without any sacrifice. On the other hand, the usage of pointbased constructions in infinite-dimensional settings requires certain ``sequential normal compactness" assumptions that are automatic in finite dimensions; see \cite{m06} for all the details. Due to this, we concentrate in what follows on neighborhood conditions, which are sufficient for our applications here without additional requirements. 

The main characterization employed below is a neighborhood one, and it addresses to  the $\al$-covering property of $\Phi\colon X\tto Y$ around the given point 
$(\ox,\oy)\in\gph\Phi$, where
\begin{equation*}
\gph\Phi:=\big\{(x,y)\in X\times Y\;\big|\;y\in\Phi(x)\big\}
\end{equation*}
is the graph of the multifunction $\Phi$. To construct an appropriate generalized derivative of $\Phi$, we start with the generalized normals to arbitrary sets in normed spaces. Given an nonempty subset $\Th\subset Z$ of a normed space $Z$ with the norm $\|\cdot\|$, the {\em prenormal cone} (known also as the regular or Fr\'echet normal cone) to $\Th$ at $\oz\in\Th$ is defined by
\begin{equation}\label{rnc}
\Hat N(\oz;\Th):=\Big\{z^*\in Z^*\;\Big|\;\disp\limsup_{z\st{\Th}{\to}\oz}\frac{\la z^*,z-\oz\ra}{\|z-\oz\|}\le 0\Big\},
\end{equation}
where $\la\cdot,\cdot\ra$ stands for the canonical paring between $Z$ and the dual space $Z^*$, and where the symbol $z\st{\Th}{\to}\oz$ indicates that $z\to\oz$ with $z\in\Th$. If $\Th$ is convex, the prenormal cone \eqref{rnc} agrees with the classical normal cone of convex analysis.

Let $\Phi\colon X\tto Y$ be a multifunction between normed spaces, and let $(\ox,\oy)\in\gph\Phi$. The {\em precoderivative} (or regular coderivative) of $\Phi$ at $(\ox,\oy)$ is defined via the prenormal cone \eqref{rnc} to the graph as a set-valued mapping $\Hat D^*\Phi(\ox,\oy)\colon Y^*\tto X^*$ with the values
\begin{equation}\label{pcod}
\Hat D^*\Phi(\ox,\oy)(y^*):=\big\{x^*\in X^*\;\big|\;(x^*,-y^*)\in\Hat N\big((\ox,\oy);\gph\Phi\big)\big\},\quad y^*\in Y^*.
\end{equation}
If $\Phi\colon X\to Y$ is single-valued around $\oy=\Phi(\ox)$ and 
Fr\'echet differentiable at $\ox$ with the derivative $\nabla\Phi(\ox)$, then $\Hat D^*\Phi(\ox)(y^*)=\{\nabla\Phi(\ox)^*y^*\}$ for all $y^*\in Y^*$ in terms of the {\em adjoint operator} to $\nabla\Phi(\ox)$. Using the robust regularizations (via the Painlev\'e-Kuratowski outer limit from points nearby) of \eqref{rnc} and \eqref{pcod}, we arrive at the (Mordukhovich, limiting) {\em normal cone} $N(\oz;\Th)$ and {\em coderivative} $D^*\Phi(\ox,\oy)$, which are robust and possess much better calculi in comparison with \eqref{rnc} and \eqref{pcod}; see \cite{m06,m18,rw} for detailed expositions and references in finite and infinite dimensions. We do not present here the exact definitions of $N(\oz;\Th)$ and $D^*\Phi(\ox,\oy)$ since only the precoderivative \eqref{pcod} is used in our formulations and proofs below, while we mention the corresponding counterparts of the results in terms of the limiting constructions.

To proceed, define by using \eqref{pcod} the {\em covering constant} for $\Phi$ at $(\ox,\oy)$ by
\begin{equation}\label{cov-const}
\begin{array}{ll}
\Hat\al(\Phi,\ox,\oy):=\disp\sup_{\eta>0}\inf\Big\{\|x^*\|\;\Big|&x^*\in\Hat 
D^*\Phi(x,y)(y^*),\;x\in B_X(\ox,\eta),\\
&y\in\Phi(x)\cap B_Y(\oy,\eta),\;\|y^*\|=1\Big\}.
\end{array}
\end{equation}

To formulate next the covering criterion that is crucial for applications in this paper, we need to recall one notion from the geometry of Banach spaces; see, e.g., \cite[Section~2.2]{m06} and the references therein. A Banach space $Z$ is {\em Asplund} if every convex continuous function defined on an open convex set $O$ in $Z$ is Fr\'echet differentiable on a dense subset of $O$. The class of Asplund spaces is well investigated in the geometric theory of Banach spaces and has been broadly employed in variational analysis and generalized differentiation. This class admits many beautiful characterizations. One of the most useful characterizations is that a Banach space $Z$ is Asplund if and only if every closed separable subspace of $Z$ has a separable dual. The collections of Asplund spaces is fairly largely including, in particular, any Banach space that admits an equivalent norm Fr\'echet differentiable at nonzero points, and hence every reflexive Banach space. On the other hand, there exists Asplund spaces not having even an equivalent norm G\^ateaux differentiable off the origin.

The following theorem is a combination of \cite[Theorem~1.54(i)]{m06} (necessity in general normed spaces) and \cite[Theorem~4.1]{m06} (sufficiency in Asplund spaces). Note that, although \cite[Theorem~1.54(i)]{m06} was formulated in Banach spaces, its  proof based on \cite[Theorems~1.43, 1.49, and 1.52]{m06} holds in an arbitrary normed space setting.

\begin{Theorem}{\bf(covering criterion).}\label{cov-cr} Let $\Phi\colon X\tto Y$ be a set-valued mapping between normed spaces, let $(\ox,\oy)\in\gph\Phi$, and let $\Hat\al(\Phi,\ox,\oy)$ be taken from \eqref{cov-const}. Then:

{\bf(i)} The $\al$-covering property of $\Phi$ around $(\ox,\oy)$ implies that $\al\le\Hat\al(\Phi,\ox,\oy)$.

{\bf(ii)} Assume that $X$ and $Y$ are Asplund spaces. If $0<\al<\Hat\al(\Phi,\ox,\oy)$, then $\Phi$ has the $\al$-covering property around $(\ox,\oy)$ provided that its graph is closed around this point. 
\end{Theorem}

It follows from Theorem~\ref{cov-cr} that the condition $\Hat\al(\Phi,\ox,\oy)>0$ is {\em necessary and sufficient} for the covering property of $\Phi$ around $(\ox,\oy)$ with {\em some} modulus $\al>0$. If both spaces $X$ and $Y$ are 
finite-dimensional, then the constant \eqref{cov-const} in the aforementioned results can be equivalently replaced by the {\em pointbased} one
\begin{equation}\label{cov-exact}
\al(\Phi,\ox,\oy):=\inf\big\{\|x^*\|\;\big|\;x^*\in D^*\Phi(\ox,\oy)(y^*),\;\|y^*\|=1\big\}
\end{equation}
expressed via the (limiting) coderivative $D^*\Phi(\ox,\oy)$ exactly at the point in question; see \cite[Theorem~3.6]{m93} and \cite[Theorem~3.3]{m18}, where the reader can also find precise formulas for characterizing the equivalent metric regularity and Lipschitz-like properties with computing the exact bounds of their moduli in finite-dimensional spaces. The results of this type are known as the {\em Mordukhovich criteria}; see \cite[Theorems~9.40 and 9.43]{rw}. Infinite-dimensional extensions of the pointbased characterizations of the covering and related properties are given in \cite[Theorems~4.10 and 4.18]{m06}. 

To conclude this section, let us recall yet another stability-type property of variational analysis that is used is what follows, together with its parametric version, in some results of the paper. We say that a multifunction $\Phi\colon X\tto Y$ between normed spaces is {\em calm} at $(\ox,\oy)\in\gph\Phi$ if \eqref{lip} holds with $u=\ox$.

If $f\colon X\to Y$ as a single-valued mapping defined on a metric space $X$ with the metric $\rho$ and values belonging to a normed space $Y$, then the calmness property of $f$
at $\ox$ with some modulus $\kk\ge 0$ is formulated as follows:
\begin{equation}\label{calm}
\|f(x)-f(\ox)\|\le\kk\rho(x,\ox)\;\mbox{ for all }\;x\in U,
\end{equation}
where $U$ is a neighborhood of $\ox$. In some situations (including those in Section~\ref{sec:marg}), one-sided calmness counterparts of \eqref{calm} for $f\colon X\to\R$ are useful as pleasable assumptions and conclusions. To this end, we say that $f$ is {\em calm at $\ox$ from above} if 
\begin{equation}\label{calm-above}
f(x)\ge f(\ox)-\kk\rho(x,\ox)\;\mbox{ when }\;x\in U.
\end{equation}
{\em Calmness from below} is defined as the fulfillment of \eqref{calm-above} for $-f$. Note that the calmness properties can be viewed as a certain Lipschitzian behavior of $f$ {\em at the point} $\ox$, but usually Lipschitzian behavior is associated with {\em two-point} properties. Although the calmness and its one-sided counterparts are not robust properties and do not possess verifiable characterizations, they are proved to be useful in some applications, mainly in finite dimensions; see, e.g., \cite{dr,m18,rw} and the references therein.

\section{Parameterized Coincidence Point Theorem}\label{sec:coin}\setcounter{equation}{0}

This section presents the main result of the paper on parameterized coincidence point of two set-valued mappings $F\colon X\tto Y$ and $G\colon X\times P\tto Y$ around given points $\ox\in X$, $\oy\in Y$ and a parameter $\op\in P$ satisfying
\begin{equation}\label{coin-bar}
\oy\in F(\ox)\cap G(\ox,\op).
\end{equation} 
Unless otherwise stated, the space of parameters $P$ is an arbitrary {\em topological space}.  We use the standard notion of {\em lower semicontinuity} of mapping of set-valued mappings on topological spaces, while observing that this notion corresponds to {\em inner semicontinuity} in the terminology of variational analysis; see \cite{rw} and \cite{m06,m18}. We say that the multifunction $F$ is {\em closed around} $(\ox,\oy)\in\gph F$ if there exist neighborhoods $U$ of $\ox$ and $V$ of $\oy$ such that the set $\gph F\cap(\cl U\times\cl V)$ is closed in $X\times Y$, where ``cl" indicates the closure operation. Here is the main result.

\begin{Theorem}{\bf(existence of parameterized coincidence points).}\label{main} Let the spaces $X$ and $Y$ in \eqref{cp} be Asplund, and let the following conditions be satisfied:
\begin{enumerate}\vspace*{-0.1in}
\item[\bf(A1)] The multifunction $F$ is closed around $(\ox,\oy)$.

\item[\bf(A2)] There are neighborhoods $U\subset X$ of $\ox$, $V\subset Y$ of $\oy$, and $O$ of $\op\in P$ as well as a number $\ell\ge 0$ such that the multifunction $x\mapsto G(\cdot,p)$ is Lipschitz-like on $U$ relative to $V$ for each $p\in O$ with  the uniform modulus $\ell$, while the multifunction $p\mapsto G(\ox,p)$ is lower/inner semicontinuous at $\op$.

\item[\bf(A3)] The Lipschitzian modulus $\ell$ of $G(\cdot,p)$ is chosen as $\ell<\hat\al(F,\ox,\oy)$, where $\hat\al(F,\ox,\oy)$ is the covering constant of $F$ around $(\ox,\oy)$ taken from \eqref{cov-const}.
\end{enumerate}\vspace*{-0.1in}
Then for each $\al>0$ with $\ell<\al<\hat\al(F,\ox,\oy)$ there exist a neighborhood $W\subset P$ of $\op$ and a single-valued mapping $\sigma\colon W\to X$ such that whenever $p\in W$ we have
\begin{equation}\label{coin-est}
F\big(\sigma(p)\big)\cap G\big(\sigma(p),p\big)\ne\emp\;\mbox{ and }\;\|\sigma(p)-\ox\|\le\disp\frac{{\rm dist}\big(\oy;G(\ox,p)\big)}{\al-\ell}.
\end{equation} 
\end{Theorem}
{\bf Proof}. Observe first that the $\al$-covering property of $\Phi\colon X\tto Y$ on $U$ relative to $V$ in \eqref{cov} can be equivalently rewritten as
\begin{equation}\label{cov1}
B_X(x,r)\subset U\Longrightarrow\big[B_Y\big(\Phi(x),\al r\big)\cap V\subset\Phi\big(B_X(x,r)\big)\big],
\end{equation}
while the $\al$-covering property of $\Phi$ around $(\ox,\oy)\in\gph\Phi$ is described in terms of \eqref{cov1} accordingly. In the proof below, we are going to use a consequence of the (nonparameterized) coincidence point theorem from \cite[Theorem~3.1]{aaz}, which is valid in complete metric spaces, concerning two set-valued mappings such that one of them has the covering property and the other is Lipschitz-like around the point in question. This theorem tells us that if $\Phi\colon X\tto Y$ enjoys the $\al$-covering property \eqref{cov1} with $U:=B_X(\ox,r)$ and $V:=B_Y(\oy,\al r)$, the set $\gph\Phi\cap(U\times V)$ is closed in $X\times Y$, while the mapping $\Psi\colon X\tto Y$ is Lipschitz-like on $U:=B_X(\ox,r)$ relative to $V:= B_Y(\oy,\al r)$ with modulus $\ell<\al$ satisfying the relationship
\begin{equation*}
{\rm dist}\big(\oy;\Psi(\ox)\big)<(\al-\ell)r,
\end{equation*}
then there exists a point $\sigma\in X$ for which
\begin{equation*}
\Phi(\sigma)\cap\Psi(\sigma)\ne\emp\;\mbox{ and }\;\|\sigma-\ox\|\le\disp\frac{{\rm dist}\big(\oy;\Psi(\ox)\big)}{\al-\ell}.
\end{equation*}

To proceed with the proof of our parameterized coincidence point theorem, pick any number $\al\in\big(\ell,\Hat\al(F,\ox,\oy)\big)$ and show that there exists a neighborhood $W$ of the nominal parameter $\op$ such that for any fixed parameter $p\in W$ we can apply to the mappings $\Phi:=F$ and $\Psi:=G(\cdot,p)$ the aforementioned result from \cite{aaz}.

Indeed, it follows from Theorem~\ref{cov-cr} that there is a number $r>0$ such that the multifunction $F$ enjoys the $\al$-covering property  on $B_X(\ox,r)$ relative to $B_Y(\oy,\al r)$. Using assumption (A1) by decreasing the radius $r$ if necessary, we get without loss of generality that the graph of $F$ is closed relative to $B_X(\ox,r)\times B_Y(\oy,\al r)$. Furthermore, decreasing again $r>0$ if necessary and using the uniform Lipschitz-like property of $G(\cdot,p)$ imposed in assumption (A2) ensure that $G(\cdot,p)$ is Lipschitz-like on $B_X(\ox,r)$ relative to $B_Y(\oy,\al r)$. It is easy to check that the inner semicontinuity assumption on the set-valued mapping $G(\ox,\cdot)$ in (A2) yields the upper semicontinuity of the parametric distance function $p\mapsto{\rm dist}(\oy;G(\ox,p))$ as $p\in W$ at $\op$. Since ${\rm dist}(\oy;G(\ox,\op))=0$ by \eqref{coin-bar}, there exists a neighborhood $O$ of $\op$ such that
\begin{equation*}
{\rm dist}\big(\oy;G(\ox,p)\big)<(\al-\ell)r\;\mbox{ for all }\;p\in O.
\end{equation*}
Therefore, for any parameter $p\in O$ the multifunctions $\Phi=F$ and $\Psi=G(\cdot,p)$ satisfy the assumptions of the aforementioned result from \cite{aaz}. This tells us that for each $p\in O$ there exists $\sigma(p)\in X$ satisfying the conditions in \eqref{coin-est}, which thus completes the proof of our main theorem. $\h$

Note that the coincidence point function $\sigma(p)$ in Theorem~\ref{main} cannot be generally selected as a {\em continuous} function with respect to the parameter $p$. This is demonstrated in the next section in the case of a single-valued continuous mapping $F\colon\R^2\to\R^2$ with the simplest parameter-dependent mapping $G(p)
\equiv p$ as $p\in\R^2$ in \eqref{pge}. 

\section{Parameterized Generalized Equations and Implicit Functions}\label{sec:pge} \setcounter{equation}{0}

In this section, we derive from Theorem~\ref{main} two results of their own interest and importance. The first result addresses standard parameterized generalized equations of type \eqref{pge} with single-valued parameter-dependent mappings $g(x,p)$. The second result gives a new implicit function theorem under general conditions. 

We start with the standard class of parameterized generalized equations where the parameter-dependent mapping in \eqref{pge} is {\em single-valued}. To relate the existence of solutions for such systems to the parameterized coincidence point problem \eqref{cp}, it is convenient to rewrite the corresponding  generalized equation \eqref{pge} in the form
\begin{equation}\label{pge1}
0\in F(x)-g(x,p),
\end{equation}
where $F\colon X\tto Y$ and $g\colon X\times P\to Y$. The next theorem is a specification of Theorem~\ref{main} for the case of \eqref{pge1} with some important addition.

\begin{Theorem}{\bf(existence of solutions to parameterized generalized equations).}
\label{A2} Let the spaces $X$ and $Y$ be Asplund, and let the triple $(\ox,\oy,\op)\in X\times Y\times P$ with $\oy:=g(\ox,\op)$ satisfy the generalized equation \eqref{pge1}.
In addition to {\rm(A1)} and {\rm(A3)}, suppose that assumption {\rm(A2)} is replaced by the following:\vspace*{-0.1in}
\begin{enumerate}
\item[\bf(A2$^{\prime}$)] There are neighborhoods $U\subset X$ of $\ox$  and $O$ of $\op\in P$ as well as a positive number $\ell$ such that the mapping $g(\cdot,p)$ is Lipschitz continuous on $U$ for each $p\in O$ with the uniform modulus $\ell$ 
while the mapping $g(\ox,\cdot)$ is continuous on $O$.
\end{enumerate}\vspace*{-0.1in}
Then for each $\al>0$ with $\ell<\al<\hat\al(F,\ox,\oy)$ there exist a neighborhood $W\subset P$ of $\op$ and a single-valued mapping $\sigma\colon W\to X$ such that whenever $p\in W$ we have
\begin{equation}\label{coin-est1}
0\in F\big(\sigma(p)\big)-g\big(\sigma(p),p\big)\;\mbox{ and }\;\|\sigma(p)-\ox\|\le\disp\frac{{\rm dist}\big(g(\ox,p);F(\ox)\big)}{\al-\ell}.
\end{equation} 
If furthermore the space $P$ is metric and the mapping $g(x,p)$ is continuous around  $(\ox,\op)$, then {\rm(A2$^{\prime}$)} is satisfied for any $\ell>\inf_{r>0}\vartheta(r)<\infty$ with $\vt(\cdot)$ given by
\begin{equation*}
\vt(r):=\sup\big\{\beta\big(g(\cdot,p),x\big)\;\big|\;x\in B_X(\ox,r),\;r\in B_P(\op,r)\big\},\quad r>0,
\end{equation*}
where the quantity $\bb(h,x)$ is defined for any $h\colon X\to Y$ and $x\in X$ via the norm of the precoderivative \eqref{pcod} as a positive homogeneous mapping by
\begin{equation}\label{cod-lip}
\bb(h,x):=\disp\inf_{t>0}\sup\big\{\|\hat D^*h(u)\|\;\big|\;u\in B_X(x,t)\big\}. 
\end{equation}
\end{Theorem}\vspace*{-0.05in}
{\bf Proof}. To verify \eqref{coin-est1}, take an arbitrary positive constant $\al$ with $\ell<\al<\hat\al(F,\ox,\oy)$, define the mapping $G\colon X\times P\tto Y$ by
\begin{equation}\label{G}
G(x,p):=\big\{-g(x,p)\}\;\mbox{ for all }\;x\in X\;\mbox{ and }\;p\in P,
\end{equation}
and observe that assumption (A2$^{\prime}$) on $g$ yields the fulfillment of (A2) for $G$. Theorem~\ref{main} gives us a neighborhood $O$ of $\op$ and a mapping $\sigma\colon O\to X$ satisfying the conditions in \eqref{coin-est} for $G$ from \eqref{G}. Taking into account the structure of \eqref{pge1}, we clearly arrive at the claimed relationships in \eqref{coin-est1}.

To prove now the second part of the theorem, use the assumption that $\ell>\inf_{r>0}\vartheta(r)$, where the latter infimum is a finite number, and then find $\ell_1,\ell_2>0$ such that 
\begin{equation*}
\ell>\ell_1>\ell_2>\disp\inf_{r>0}\vt(r).
\end{equation*}
The definition of $\vt(r)$ ensures the existence of $r>0$ for which
\begin{equation*}
\bb\big(g(\cdot,p),x\big)\le\ell_2\;\mbox{ whenever }\;p\in B_P(\op,r)\;\mbox{ and }\;x\in B_X(\ox,r).
\end{equation*}
In our further discussions, suppose that $U$ is the interior of the ball $B_X(\ox,r)$ and $O$ is the interior of the ball $B_P(\op,r)$. Fix $p\in O$ and pick any $x\in U$. Using the assumption $\bb(g(\cdot,p),x)\le\ell_2$ and employing the neighborhood characterization of the local Lipschitz continuity taken from \cite[Theorem~4.7]{m06} (which is actually equivalent of Theorem~\ref{cov-cr} for inverse mappings), we conclude by \eqref{cod-lip} that the mapping $g(\cdot,p)$ is locally Lipschitzian around $x$ with modulus $\ell_1>\ell_2$. Since the latter holds for any $x\in U$ with the same constant $\ell_1$, this verifies that $g(\cdot,p)$ is Lipschitz continuous on $U$ with the given modulus $\ell$. This verifies, due to the arbitrary choice  of $p\in O$, the fulfillment of assumption (A2) and thus completes the proof of the theorem. $\h$ 

The following example shows that the solution $\sigma(p)$ to \eqref{pge1}, and hence to the coincidence point problem in \eqref{cp}, cannot be generally selected as a {\em continuous} function even in the finite-dimensional setting with single-valued smooth mappings $F$ and $g$. 

\begin{Example}{\bf(discontinuity of solutions to parameterized generalized equations).}\label{discont} There exist smooth mappings $F\colon\R^2\to\R^2$ and $g\colon\R^2\times\R^2\to\R^2$ satisfying all the assumptions of Theorem~{\rm\ref{A2}} such that any local solution $\sigma(p)$ to the parameterized generalized equation \eqref{pge1} having properties \eqref{coin-est1} is discontinuous.
\end{Example}\vspace*{-0.05in}
{\bf Proof}. Define the smooth mapping $F\colon\R^2\to\R^2$ by 
\begin{equation*}
F(x):=\disp\frac{1}{2}\left(\begin{array}{ll}x^2_1-x^2_2\\
\;\;2x_1 x_2
\end{array}\right)
\;\mbox{ for }\;x=(x_1,x_2)\ne 0\;\mbox{ with }\;F(0)=0.
\end{equation*}
It is shown in \cite{az} that $\Hat\al(F,0,0)=1/2$ for the covering constant 
\eqref{cov-const}, while any (right) inverse mapping $\sigma(\cdot)$ to $F$ near $0\in\R^2$ with $\sigma(0)=0$ is discontinuous. Letting $g(x,p)\equiv p$, $\ell=0$, and $\ox=\op=0$ tells us that all the assumptions of Theorem~\ref{A2} are 
satisfied. The corresponding mapping $\sigma(p)$ is inverse to $F$ with the condition $\sigma(0)=0$. This shows that the continuity of solutions fails in this problem. $\h$

The next application of Theorem~\ref{main} concerns {\em implicit function theorems}, which--together with related {\em inverse function} results--constitute a fundamental topic in analysis and applications, including their variational aspects; see, e.g., the book \cite{dr} with the references therein. In the framework of this paper, we formulate the implicit function problem as follows. Given a mapping $f
\colon X\times P\to Y$, consider the equation
\begin{equation}\label{impl}
f(x,p)=0
\end{equation}
with the decision variable $x\in X$ and the parametric variable $p$ belonging to a topological space $P$. Fixing a nominal parameter value $p=\op$, suppose that a point $\ox$ solves equation \eqref{impl} when $p=\op$, i.e.,
\begin{equation}\label{imp0}
f(\ox,\op)=0.
\end{equation}
We are interested  in finding efficient conditions that ensure the existence of {\em single-valued} solutions $\sigma(p)$ to \eqref{impl} for all parameters $p$ near $\op$ with some constructive {\em error bound}, provided that the mapping $f(\cdot,\op)$ enjoys the {\em covering property} around $\ox$. Note that the coderivative characterizations of covering have been used in variational analysis to establish the existence of {\em implicit multifunctions} exhibiting a certain local Lipschitzian stability; see, in particular, \cite[Corollaries~4.38 and 4.42]{m06}. However, we are not familiar with the usage of covering constants to get {\em single-valued implicit functions}. In fact, this was one of our initial motivations for this paper.

Here is the implicit function theorem, which follows from Theorem~\ref{main}.

\begin{Theorem}{\bf(implicit function theorem from covering).}\label{impl-cov} Let the pair $(\ox,\op)$ satisfy \eqref{imp0}, where the spaces $X$ and $Y$ are Asplund. Suppose that the mapping $f(\cdot,\op)$ is closed around $(\ox,f(\ox,\op))$, while the mapping $f(\ox,\cdot)$ is continuous at $\op$. Assume in addition that there are neighborhoods $U\subset X$ of $\ox$ and $O\subset P$ of $\op$ such that
\begin{equation}\label{lip-est1}
\big\|\big(f(x_1,p)-f(x_1,\op)\big)-\big(f(x_2,p)-f(x_2,\op)\big)\big\|\le\ell\|x_1-x_2\|
\end{equation}
for all $x_1,x_2\in U$ and $p\in O$. If $f(\cdot,\op)$ has the covering property around $\ox$, i.e.,
\begin{equation*}
\Hat\al\big(f(\cdot,\op),\ox\big)>0
\end{equation*}
for the covering constant from \eqref{cov-const}, then there exist a neighborhood $W$ of $\op$, a number $c>0$, and a mapping $\sigma\colon W\to X$ such that
\begin{equation}\label{impl-est}
f\big(\sigma(p),p\big)=0\;\mbox{ and }\;\|\sigma(p)-\ox\|\le c\|f(\ox,p)\|\;\mbox{ for all }\;p\in W.
\end{equation}
\end{Theorem}\vspace*{-0.1in}
{\bf Proof}. Pick any numbers $\al>0$ and $\ell\ge 0$ with $\al<\Hat\al(f(\cdot,\op),\ox)$ and $\ell<\al$. Define
\begin{equation*}
F(x):=\big\{f(x,\op)\big\}\;\mbox{ and }\;G(x,p):=\big\{g(x,p)-g(x,\op)\big\}\;\mbox{ as }\;x\in X,\;p\in P
\end{equation*}
and show that these mappings satisfy the assumptions of Theorem~\ref{main}. Indeed, (A1) holds since $f(\cdot,\op)$ is closed around $(\ox,f(\ox,\op))$. It is easy to check that \eqref{lip-est1} yields (A2). Assumption (A3) follows from the construction of $F$ and $G$ due to the covering characterization of Theorem~\ref{cov-cr}. Applying now Theorem~\ref{main}, we get \eqref{coin-est}, which clearly yields \eqref{impl-est} with $c:=(\al-\ell)^{-1}$ in the setting under consideration. $\h$

\section{Optimal Value Functions}\label{sec:marg}
\setcounter{equation}{0}

This section is devoted to studying a broad class of optimal value functions associated with parametric minimization problems subject to parameterized generalized equation constraints. Given a set-valued mapping $F\colon X\tto Y$, single-valued mapping $g\colon X\times P\to Y$, and an extended-real-valued function $\ph\colon X\times P\to\oR$, consider the following family of parametric optimization  problems with parameterized generalized equation constraints written in form \eqref{pge1}:
\begin{equation}\label{po}
\mbox{minimize }\;\ph(x,p)\;\mbox{ subject to }\;0\in F(x)-g(x,p)
\end{equation}
depending on the parameter $p$ from a topological space $P$. The {\em optimal value function} for the parametric problem \eqref{po} is defined by
\begin{equation}\label{marg}
\mu(p):=\inf\big\{\ph(x,p)\;\big|\;0\in F(x)-g(x,p)\big\},\quad p\in P.
\end{equation}
It has been well recognized in variational analysis and optimization that optimal value functions (known also as {\em marginal functions}) play a highly important role in various qualitative and numerical aspects of the theory and applications; see, e.g., \cite{m06,m18,rw} with the extensive commentaries and bibliographies therein. However, not much has been done in the framework of parameterized generalized equations as in \eqref{marg}. 

The main goal here is investigate the (lower and upper) {\em semicontinuity} and {\em continuity} properties, as well as some Lipschitzian counterparts, of optimal value functions of type \eqref{marg}. These issues are undoubtedly crucial to conduct, e.g., {\em sensitivity} and {\em stability analysis} of optimization problems under parameter perturbations. Our approach is based on the solvability result for the parameterized generalized equations with the error bound estimate obtained in Theorem~\ref{A2}, which in turn utilizes the covering characterization taken from Theorem~\ref{cov-cr}. Note that, although we cannot guarantee the continuity of parameterized solutions in Theorem~\ref{A2}, the results obtained in that theorem occur to be instrumental to derive efficient conditions for the semicontinuity and continuity of the optimal value function \eqref{marg}. 

The major results of this section separately establish conditions for the {\em lower semicontinuity} and {\em upper semicontinuity} of the optimal value function \eqref{marg} with the corresponding assertions on its {\em calmness from below} and {\em calmness from above}. Thus the unification of these results justifies the {\em continuity} and the {\em calmness} of \eqref{marg}. We also find efficient conditions ensuring the local {\em Lipschitz continuity} of the optimal value function around the reference point. 

Let $\ox$ solve problem \eqref{po} with the fixed parameter value $p=\op$, i.e.,
\begin{equation}\label{po1}
\mbox{minimize }\;\ph(x,\op)\;\mbox{ subject to }\;0\in  F(x)-g(x,\op)
\end{equation}
meaning that $0\in F(\ox)-g(\ox,\op)$ and $\ph(\ox)\le\ph(x)$ if $0\in F(x)-g(x,\op)$. This yields
\begin{equation*}
\mu(\op)=\ph(\ox)>-\infty.
\end{equation*}
The first theorem of this section provides conditions under which the optimal value function is {\em upper semicontinuous} at $\ox$ and {\em calm from above} at this point.

\begin{Theorem}{\bf(upper semicontinuity and calmness from above of optimal value functions).}\label{usc} Let $\ox$ be a solution to problem \eqref{po1}, where $\ph$ is upper semicontinuous at $(\ox,\op)$. Impose the following assumptions:\vspace*{-0.1in}
\begin{itemize}
\item The multifunction $F\colon X\tto Y$ between Asplund spaces is closed around $(\ox,0)$.
\item There are neighborhoods $U$ of $\ox$ and $O$ of $\op$ as well as a number $\ell\ge 0$ such that the mapping $g(\cdot,p)$ is Lipschitz continuous on $U$ for each $p\in O$ with the uniform modulus $\ell$, while the mapping $g(\ox,\cdot)$ is continuous on $O$.
\item We have $\ell<\Hat\al\big(F,\ox,g(\ox,\op)\big)$ for the covering constant from \eqref{cov-const}.
\end{itemize}\vspace*{-0.1in}
Then the optimal value function $\mu(\cdot)$ from \eqref{marg} is upper semicontinuous at $\op$.

If furthermore $P$ is a metric space with the metric $\rho$, and if the cost function $\ph$ and the mapping $g(\ox,\cdot)$ are calm at the points $(\ox,\op)$ and $\op$, respectively, with some moduli in \eqref{calm}, then there exists $\kk>0$ such that
\begin{equation}\label{calm1}
\mu(p)\le\mu(\op)+\kk\rho(p,\op),
\end{equation}
i.e., the optimal value function \eqref{marg} is calm from above at $\op$ in the sense of \eqref{calm-above}. 
\end{Theorem}\vspace*{-0.05in}
{\bf Proof}. Pick arbitrary numbers $\al>0$ and $\ell\ge 0$ together with neighborhoods $U$ of $\ox$ and $O$ of $\op$ such that $\ell<\al<\Hat\al\big(F,\ox,-g(\ox)\big)$, the mapping $g(\cdot,p)$ is Lipschitz continuous on $U$ for all $p\in O$ with the uniform modulus $\ell$, and $g(\ox,\cdot)$ is continuous on $U$. We see that all the assumptions in the first part of Theorem~\ref{A2} are satisfied, and thus there exists a neighborhood $W$ of $\op$ and a mapping $\sigma\colon W\to X$  for which both relationships in \eqref{coin-est1} hold. The first condition therein tells us that the point $\sigma(p)$ is feasible in \eqref{po} for each $p\in W$. Therefore, we have
\begin{equation}\label{feas}
\mu(p)\le\ph\big(\sigma(p),p\big)\;\mbox{ whenever }\;p\in W.
\end{equation}

Let us prove that $\mu(\cdot)$ is upper semicontinuous at $\op$. Fix any $\ve>0$ and using the assumed upper semicontinuity of $\ph$ at $(\ox,\op)$, find a number $\dd>0$ and a neighborhood $O_1\subset W$ of the point $\op$ such that 
\begin{equation}\label{usc1}
\ph(x,p)<\ph(\ox,\op)+\ve\;\mbox{ for all }\;x\in B_X(\ox,\dd)\;\mbox{ and }\;p\in O_1.
\end{equation} 
The continuity assumption on $g(\ox,\cdot)$ allows us to get a neighborhood 
$O_2\subset O_1$ with
\begin{equation*}
{\rm dist}\big(g(\ox,p);F(\ox)\big)\le(\al-\ell)\dd\;\mbox{ for all }\;p\in O_2.
\end{equation*}
Employing this together with the estimate in \eqref{coin-est1} gives us
\begin{equation*}
\|\sigma(p)\|\le\disp\frac{{\rm dist}\big(g(\ox,p);F(\ox)\big)}{\al-\ell}\le\dd 
\end{equation*}
and thus verifies the inclusion
\begin{equation}\label{usc2}
\sigma(p)\in B_X(\ox,\dd)\;\mbox{ whenever }\;p\in O_2.
\end{equation}
Combining \eqref{usc1} and \eqref{usc2}, we arrive at the relationships
\begin{equation*}
\mu(p)\le\ph\big(\sigma(p),p\big)<\ph(\ox,\op)+\ve=\mu(\op)+\ve\;\mbox{ for all }\;p\in O_2,
\end{equation*} 
which justify the claimed upper semicontinuity of $\mu(\cdot)$ at $\op$. 

Now we proceed with the verifying that the optimal value function \eqref{marg} enjoys the calmness from above property \eqref{calm1} at $\op$ under the additional assumptions imposed in the theorem. These assumptions allow us to find $\dd>0$ and a neighborhood $O_2\subset W$ of $\op$ ensuring the inclusion \eqref{usc2} together with the estimate
\begin{equation}\label{usc3}
|\ph(x,p)-\ph(\ox,\op)|\le\kk_1\big(\|x-\ox\|+\rho(p,\op)\big)\;\mbox{ for all }\;p\in O_2,\;\;x\in B_X(\ox,r)
\end{equation}
with some constant $\kk_1>0$. Taking into account that $0\in F(\ox)-g(\ox,\op)$ leads us to
\begin{equation}\label{usc4}
{\rm dist}\big(g(\ox,p);F(\ox)\big)\le\|g(\ox,p)-g(\ox,\op)\|\le{\rm dist}\big(g(\ox,\oy);F(\ox)\big)\le\kk\rho(p,\op)
\end{equation}
whenever $p\in O_2$. Thus for all such $p$ we get the estimates
\begin{equation*}
\begin{array}{ll}
\mu(p)&\le\ph\big(\sigma(p),p\big)\le\ph(\ox,\op)+\kk_1\|\ox-\sigma(p)\|+\kk\rho(p,\op)\\\\
&\le\ph(\ox,\op)+\kk_1\disp\frac{{\rm dist}\big((\ox,p);F(\ox)\big)}{\al-\ell}+\kk_1\rho(p,\op)\\\\
&\le\ph(\ox,\op)+\disp\frac{\kk_1\kk_2\rho(\op,p)}{\al-\ell}+\kk_1\rho(p,\op)\le\mu(\op)+\kk\rho(p,\op),
\end{array}
\end{equation*}
where $\kk:=\kk_1+\kk_1\kk_2(\al-\ell)^{-1}$. Indeed, the first claimed inequality follows from definition \eqref{marg}, the second one is a consequence of \eqref{usc2} and \eqref{usc3}, the third and forth inequalities follow from \eqref{coin-est1} and \eqref{usc4}, respectively, while the last one is due to $\mu(\op)=\ph(\ox)$ and the definition of $\kk$. This verifies that $\mu(\cdot)$ is calm from above as in \eqref{calm1} and therefore completes the proof of the theorem. $\h$\vspace*{0.05in}

Next we find efficient conditions ensuring the {\em lower semicontinuity} and also {\em calmness from below} of the optimal value function \eqref{marg} at the point in question. The following simple one-dimensional example shows that the lower semicontinuity of $\mu(\cdot)$ may fail under the assumptions of Theorem~\ref{usc}.

\begin{Example}{\bf(failure of lower semicontinuity of optimal value functions).}\label{no-lsc} {\rm Let $X=Y=P:=\R$, $F(x)\equiv\R$, $g(x,p)\equiv 0$, and
\begin{equation*}
\ph(x,p):=\max\big\{-1,|p| x\big\}\;\mbox{ for all }\;x,p\in\R.
\end{equation*}
Then all the assumptions of Theorem~\ref{usc} hold with $\ox=\op=0$. It is easy to see that $\mu(p)=-1$ for all $p\ne 0$, while $\mu(0)=0$. Thus the optimal value function \eqref{marg} is not lower semicontinuous at $\op=0$.}
\end{Example}

Observe that the lower semicontinuity is a {\em major property} of extended-real-valued functions in variational analysis and optimization inclined towards minimization. A crucial role of this property for optimal value (marginal) functions has been revealed in the books \cite{m06,m18,rw} and the references therein  dealing with problems in finite-dimensional and Banach spaces. 

To investigate now the lower semicontinuity together with the calmness from below properties of $\mu(\cdot)$ in the topological and metric space settings, consider the parametric set of feasible solutions to problem \eqref{po} given by
\begin{equation}\label{S}
S(p):=\big\{x\in X\;\big|\;0\in F(x)-g(x,p)\big\},\quad p\in P,
\end{equation}
and introduce the {\em semilocal covering constant} of $F$ depending only on the reference point $\ox\in X$ uniformly in $\oy\in F(\ox)$ as
\begin{equation*}
\Hat\al(F,\ox):=\disp\sup_{\eta>0}\inf\big\{\|x^*\|\;\big|\;x^*\in\Hat D^*F(x,y)(y^*),\;x\in B_X(\ox,\eta),\;y\in F(x),\;\|y^*\|=1\big\}.
\end{equation*}
We obviously have that $\Hat\al(F,\ox)\le\Hat\al(F,\ox,\oy)$ for any $\ox\in X$ and $\oy\in F(\ox)$. Here is the main result of this section providing sufficient conditions for the lower semicontinuity of optimal value function \eqref{marg} in general framework. 

\begin{Theorem}{\bf(lower semicontinuity and calmness from below of optimal value functions).}\label{lsc} Let $F\colon X\tto Y$ be a closed-graph mapping between Asplund spaces, let $S(\op)\ne\emp$ for the feasible solution map \eqref{S} at the nominal parameter value $\op\in P$, and let there exist numbers $\al>0$, $\ell\ge 0$, and $\dd>0$ together with a neighborhood $O$ of $\op$ such that the following assumptions are satisfied:\vspace*{-0.1in}
\begin{itemize}
\item We have $\inf\big\{\Hat\al(F,\ox)\;\big|\;x\in B_X(u,\dd),\;u\in S(\op)\big\}>\al>\ell$.
\item The mapping $g(\cdot,\op)$ is Lipschitz continuous on $B_X(u,\dd)$ with the uniform modulus $\ell$ for each $u\in S(p)$ and $p\in O$.
\item The mapping $g(x,\cdot)$ is continuous at $\op$ uniformly in $x\in X$.
\item The cost function $\ph$ is uniformly lower semicontinuous at $(x,\op)$ with $x\in S(\op)$, i.e., for any $\ve>0$ there are $\nu>0$ and a neighborhood $O$ of $\op$ such that 
\begin{equation*}
\ph(x,\op)>\ph(u,p)-\ve\;\mbox{ for all }\;x\in S(\op),\;p\in O,\;u\in X\;\mbox{ with }\;\|u-x\|\le\dd.
\end{equation*}
\end{itemize}\vspace*{-0.1in}
Then the optimal value function $\mu(\cdot)$ is lower semicontinuous at $\op$. 

Impose furthermore the following additional assumptions:\vspace*{-0.1in}
\begin{itemize}
\item $P$ is a metric space with the metric $\rho$.
\item The cost function $\ph$ is uniformly calm  at $(x,\op)$ for all $x\in S(\op)$ with respect to both variables, i.e., there exists a number $\kk_1\ge 0$ such that
\begin{equation*}
|\ph(u,p)-\ph(x,\op)|\le\kk_1\|u-x\|+\kk_1\rho(p,\op)\;\mbox{ for all }\;u\in B_X(x,\dd),\;p\in O.
\end{equation*}
\item The mapping $g(u,\cdot)$ is calm at $\op$ uniformly in $u\in B_X\big(S(\op),\dd)\big)$, i.e., there exists a number $\kk_2\ge 0$ such that
\begin{equation*}
\|g(u,p)-g(u,\op)\|\le\kk_2\rho(p,\op)\;\mbox{ for all }\;u\in B_X(\ox.\dd),\;x\in S(\op).
\end{equation*}
\end{itemize}\vspace*{0.05in}
Then the optimal value function $\mu(\cdot)$ is calm from below at $\op$, i.e., there exists a constant $\kk>0$ such that
\begin{equation}\label{calm below}
\mu(p)\ge\mu(\op)-\kk\rho(p,\op)
\end{equation}
for all $p\in O$ sufficiently close to $\op$.
\end{Theorem}\vspace*{-0.1in}
{\bf Proof}. First we verify the lower semicontinuity of $\mu(\cdot)$ under the assumptions made. The claimed lower semicontinuity means that for any $\ve>0$ there exists a neighborhood $W$ of $\op$ on which $\mu(\op)>\mu(p)-\ve$. To proceed, we employ the assumed lower semicontinuity of the cost function $\ph$ in $x$ uniformly in $p$ and find a number $\nu>0$ and a neighborhood $O\subset P$ of $\op$ such that
\begin{equation}\label{lsc1}
\ph(x,\op)>\ph(u,p)-\ve\;\mbox{ for all }\;x\in S(\op),\;u\in X\;\mbox{ with }\;\|u-x\|\le\nu,\;p\in O.
\end{equation} 
Using further the continuity of $g(x,\cdot)$ at $\op$ uniformly in $x\in X$ gives us a neighborhood $W\subset O$ of $\op$ on which we have
\begin{eqnarray}\label{lsc2}
\|g(u,p)-g(u,\op)\|\le(\al-\ell)\min\big\{\dd,\nu\big\}\;\mbox{ for any }\;p\in W\;\mbox{ and }\;u\in S(p).
\end{eqnarray}
Fix now arbitrary points $p\in W$ and $u\in S(p)$, and then show that the multifunction $F$ and the mapping $g(\cdot,\op)$ satisfy the assumptions of the aforementioned consequence of the coincidence point theorem from \cite[Theorem~3.1]{aaz} relative to the balls $B_X(u,\dd)$ and $B_Y(g(u,\op),\al\dd)$. Indeed, we have 
\begin{equation*}
\inf\big\{\Hat\al(F,x)\;\big|\;x\in B_X(u,\dd)\big\}>\al,
\end{equation*}
and hence $F$ enjoys, by Theorem~\ref{cov-cr}, the $\al$-covering property around any $x\in B_X(u,\dd)$. This tells us that $F$ possesses the latter property on $B_X(u,\dd)$ relative to $B_Y(g(u,\op),\al\dd)$. Furthermore, the Lipschitz continuity of $g(\cdot,\op)$ allows us to get 
\begin{equation}\label{lsc3}
\begin{array}{ll}
{\rm dist}\big(g(u,\op);F(x)\big)&\le{\rm dist}\big(g(u,p);F(x)\big)+\|g(u,p)-g(x,\op)\|\\
&=\|g(u,p)-g(x,\op)\|,
\end{array}
\end{equation}
where the inequality follows from the triangle inequality for the distance function on $Y$, while the equality is a consequence of the inclusion $u\in S(p)$, i.e., $g(u,p)\in F(u)$. Thus we deduce from \eqref{lsc2} the estimate
\begin{equation*}
{\rm dist}\big(g(u,\op);F(u)\big)<(\al-\ell)\dd. 
\end{equation*}
Therefore, $F$ and $g(\cdot,\op)$ satisfy the assumptions of the coincidence point theorem from \cite{aaz} relative to $B_X(u,\dd)$ and $B_Y(g(u,\op),\al\dd)$, which gives us $\sigma=\sigma(u,\op)$ such that 
\begin{equation}\label{lsc4}
0\in F(\sigma)-g(\sigma,\op)\;\mbox{ and }\;\|u-\sigma\|\le\disp\frac{{\rm dist}\big(g(u,\op);F(u)\big)}{\al-\ell}. 
\end{equation}
Thus we get that $\sigma\in S(\op)$ with the fulfillment of the estimate
\begin{equation*}
\|u-\sigma\|\le\disp\frac{{\rm dist}\big(g(u,\op);F(u)\big)}{\al-\ell}\le\disp\frac{\|g(u,p)-g(u,\op)\|}{\al-\ell}<\nu,
\end{equation*}
where the second inequality follows from \eqref{lsc3}, while the third one is a consequence of \eqref{lsc2}. Employing finally \eqref{lsc1} yields $\ph(\sigma,\op)>\ph(u,p)-\ve$, which actually holds by the proof above for any $p\in W$ and $u\in S(p)$ with $\sigma=\sigma(u,p)\in S(p)$. This readily implies that $\mu(\op)>\mu(p)-\ve$ for the optimal value function \eqref{marg} meaning that it is lower semicontinuous at the nominal parameter point $\op$. 

To complete the proof of the theorem, it remains to verify that the optimal value function is calm from below as in \eqref{calm} under the additional assumptions listed in the theorem formulation. Fixing any $p\in W$ and $u\in S(p)$, take $\sigma(u,p)$ from the coincidence point theorem above and then get
\begin{equation*}
\begin{array}{ll}
\ph(u,p)\ge\ph(\sigma,\op)-|\ph(u,p)-\ph(\sigma,\op)|\ge\mu(\op)-\kk_1\|u-\sigma\|-\kk_1\rho(p,\op)\\\\
\ge\mu(\op)-\kk_1\disp\frac{{\rm dist}\big(g(u,\op);F(u)\big)}{\al-\ell}-\kk_1\rho(p,\op)\ge\mu(\op)-\kk_1\disp\frac{\|g(u,p)-g(u,\op)\|}{\al-\ell}-\kk_1\rho(p,\op)\\\\
\ge\mu(\op)-\kk_1\disp\frac{\kk_2\rho(p,\op)}{\al-\ell}-\kk_1\rho(p,\op)=\mu(\op)-\kk\rho(p,\op),
\end{array}
\end{equation*}
where $\kk:=\kk_1+\kk_1\kk_2(\al-\ell)^{-1}$, the second inequality follows from the definition of $\mu(\cdot)$ and the uniform calmness from below of $\ph$ at $(x,\op)$ as $x\in S(\op)$, the third and forth ones are consequences of \eqref{lsc4} and \eqref{lsc3}, respectively, the fifth inequality is due to the calmness property of $g(u,\cdot)$  uniformly in $u\in B_X(S(\op,\dd)$, while the equality holds due to the above definition of $\kk$.  Since $u\in S(p)$ was chosen arbitrarily, we get \eqref{calm below} by the construction of $\mu(\cdot)$ in \eqref{marg} and thus finish the proof of the theorem. $\h$\vspace*{0.05in}

As direct consequences of Theorems~\ref{usc} and \ref{lsc}, we formulate the following two results ensuring, respectively, the continuity and calmness of the optimal value function \eqref{marg}. Note that the specification of the calmness constant in Corollary~\ref{calmness} follows from the proofs of these theorems.

\begin{Corollary}{\bf(continuity of optimal value functions).}\label{cont} Let $F\colon X\tto Y$ be a closed mapping between Asplund spaces, let $\ox\in S(\op)$ be a solution to problem \eqref{po1}, and let there exist numbers $\al>0$, $\ell\ge 0$, $\dd>0$ and a neighborhood $O$ of $\op$ such the following conditions are satisfied:\vspace*{-0.1in}
\begin{itemize}
\item $\inf\big\{\Hat\al(F,x)\;\big|\;x\in B_X(u,\dd),\;u\in S(\op)\big\}>\al>\ell$.
\item The mapping $g(\cdot,\op)$ is Lipschitz continuous at $\op$ uniformly in $x\in X$.
\item The cost function $\ph$ is uniformly lower semicontinuous at $(x,\op)$ with $x\in S(\op)$.
\end{itemize}\vspace*{-0.1in}
Then the optimal value function $ \mu(\cdot)$ is continuous at $\op$. 
\end{Corollary}

\begin{Corollary}{\bf(calmness of optimal value functions).}\label{calmness} In addition to the assumptions of Corollary~{\rm\ref{cont}}, suppose that:\vspace*{-0.1in}
\begin{itemize}

\item $P$ is a metric space with the metric $\rho$.

\item The cost function $\ph$ is uniformly calm  at $(x,\op)$ for all $x\in S(\op)$ with respect to both variables with some modulus $\kk_1\ge 0$.

\item The mapping $g(u,\cdot)$ is calm at $\op$ uniformly in $u\in B_X\big(S(\op),\dd)\big)$ with some modulus $\kk_2\ge 0$. 
\end{itemize}\vspace*{-0.1in}
Then the optimal value function is calm at $\op$ with the modulus $\kk:=\kk_1\kk_1\kk_2(\al-\ell)^{-1}$, i.e., there exists a neighborhood $O$ of $\op$ such that
\begin{equation}\label{marg-calm}
|\mu(p)-\mu(\op)|\le\kk\rho(p,\op)\;\mbox{ for all }\;p\in O.
\end{equation}
\end{Corollary}\vspace*{0.05in}

The final result of this paper presents sufficient conditions ensuring the {\em local Lipschitz continuity} of the optimal value function, i.e., the existence of a number $\kk\ge 0$ and a neighborhood $O$ of $\op$ such that a counterpart of \eqref{marg-calm} holds with the replacement of $\op$ by any point $s\in O$. To furnish this, we introduce the following convexity property in the setting of normed spaces $X,Y,P$ in which case $\rho$ is the norm in $P$. 

Given mappings $F\colon X\tto Y$ and $g\colon X\times P\to Y$ between normed spaces, we say that the {\em pair} $(g,F)$ is {\em convex} if 
\begin{equation}\label{convex}
\lm g(x_1,p_1)+(1-\lm)g(x_2,p_2)\in F\big(\lm x_1+(1-\lm)x_2\big)
\end{equation}
whenever $g(x_i,p_i)\in F(x_i)$ as $i=1,2$ and $\lm\in[0,1]$. Note that condition \eqref{convex} holds for {\em any} mapping $g$ provided that 
\begin{equation}\label{convex1}
F(x+z)\supset F(x)+F(z)\;\mbox{ and }\;F(\lm x)\supset\lm F(x)\;\mbox{ for any }\;x,z\in X,\lm\ge 0.
\end{equation}
If the second inclusion in \eqref{convex1} holds as equality, then \eqref{convex1} is known as a {\em convex process} in the sense of Rockafellar; see \cite[Section~39]{r} for various properties of convex processes in finite dimensions. It is obvious that the pair convexity in \eqref{convex}, which takes into account the mapping $g$, is essentially more general and flexible for applications than merely the process convexity of $F$ as in \eqref{convex1}.   

Now we are in a position to establish conditions for local Lipschitz continuity of the class of optimal value functions in \eqref{marg}.

\begin{Theorem}{\bf(Lipschitz continuity of optimal value functions).}\label{marg-lip} In addition to the assumptions of Corollary~{\rm\ref{cont}}, suppose that the parameter space $P$ is normed, that the cost function $\ph$ is convex, and that the pair $(g,F)$ is convex in the sense of \eqref{convex}. Then the optimal value function $\mu(\cdot)$ is locally Lipschitzian around $\op$. 
\end{Theorem}\vspace*{-0.06in}
{\bf Proof}. We can obviously rewrite \eqref{marg} in the form
\begin{equation}\label{margS}
\mu(p)=\inf\big\{\ph(x,p)\;\big|\;x\in S(p)\big\},\quad p\in P,
\end{equation}
where $S$ is taken from \eqref{S} with $S(\op)\ne\emp$ by $0\in F(\ox)-g(\ox,\op)$.

First we show that $\mu(\cdot)$ in \eqref{margS} is convex provided that $\ph$ is convex together with the graph of $S$. To proceed, take any $p_1,p_2$ where $\mu(\cdot)$ is finite and fix any $\lm\in(0,1)$. Given an arbitrary positive number $\ve$, find $x_i\in S(p_i)$ as $i=1,2$ with
\begin{equation*}
\ph(x_i,p_i)>\mu(p_i)+\ve\;\mbox{ for }\;i=1,2.
\end{equation*}
This directly yields the strict inequalities
\begin{equation*}
\lm\ph(x_1,p_1)<\lm\mu(p_1)+\lm\ve\;\mbox{ and }\;(1-\lm)\ph(x_2,p_2)<(1-\lm)\mu(p_2)+(1-\lm)\ve.
\end{equation*}
Summing up these inequalities and using the convexity of $\ph$, we get
\begin{equation*}
\begin{array}{ll}
\ph\big(\lm x_1+(1-\lm)x_2,\lm p_1+(1-\lm)p_2\big)&\le\lm\ph(x_1,p_1)+(1-\lm)\ph(x_2,p_2)\\
&<\lm\mu(p_1)+(1-\lm)\mu(p_2)+\ve. 
\end{array}
\end{equation*} 
Furthermore, invoking the convexity of the graph of $S$ tells us that
\begin{equation*}
\big(\lm x_1+(1-\lm)x_2,\lm p_1+(1-\lm)p_2\big)=\lm(x_1,p_1)+(1-\lm)(x_2,p_2)\in\gph S,
\end{equation*}
and therefore we have $\lm x_1+(1-\lm)x_2\in S(\lm p_1+(1-\lm)p_2)$, which implies that
\begin{equation*}
\begin{array}{ll}
\mu\big(\lm p_1+(1-\lm)p_2\big)&\le\ph\big(\lm x_1+(1-\lm)x_2,\lm p_1+(1-\lm)p_2\big)
\\
&<\lm\mu(p_1)+(1-\lm)\mu(p_2)+\ve.
\end{array}
\end{equation*}
Since $\ve>0$ was chosen arbitrarily, this verifies the convexity of $\mu(\cdot)$ under the convexity assumptions on $\ph$ and $\gph S$.

Now we check that the graph of $S$ is convex if the pair $(g,F)$ is convex in the sense of \eqref{convex}. Observe the representations
\begin{equation}\label{Sgph}
\begin{array}{ll}
\gph S&=\big\{(p,x)\in P\times X\;\big|\;g(x,p)\in F(x)\big\}\\
&=\{(p,x)\in P\times X\;\big|\;\big(x,g(x,p)\big)\in\gph F\big\}.
\end{array}
\end{equation} 
To verify the convexity of $\gph S$, take any pairs $(p_i,x_i)\in P\times X$ such that $(p_i,x_i)\in\gph S$ as $i=1,2$. By the second representation of $\gph S$ in \eqref{Sgph}, this means that 
\begin{equation}\label{Sgph1}
\big(x_i,g(x_i,p_i)\big)\in\gph F\;\mbox{ for }\;i=1,2.
\end{equation} 
Fix any $\lm\in[0,1]$ and consider the convex combination
\begin{equation}\label{gph1}
\begin{array}{ll}
&\lm\big(x_1,g(x_1,p_1)\big)+(1-\lm)\big(x_2,g(x_2,p_2)\big)\\
&=\big(\lm x_1+(1-\lm)x_2,\lm g(x_1,p_1)+(1-\lm)g(x_2,p_2)\big). 
\end{array}
\end{equation}
Using \eqref{Sgph1} and the assumed convexity of the pair $(g,F)$ from \eqref{convex} readily ensures that the convex combination in \eqref{gph1} belongs to the graph of $F$. Therefore, the latter graph is a convex set, and so is the graph of $S$ by \eqref{Sgph}. This verifies the convexity of the optimal value function $\mu(\cdot)$ via representation \eqref{margS}. 

Finally, it has been well recognized in convex analysis that if a convex function defined on a normed space is continuous at some point, then it is Lipschitz continuous around this point; see, e.g., \cite[Corollary~2.150]{mn} and also Theorem~2.144 therein for characterizations of these equivalent properties. Unifying this with the continuity result of Corollary~\ref{cont}, we complete the proof of the theorem. $\h$

\section{Conclusions}\label{conc}

This paper shows that appropriate tools of variational analysis and generalized differentiation married to coincidence point results in nonparameterized settings lead us to establishing novel existence theorems in a broad class of parameterized  generalized equations that are highly important for various aspects of parametric optimization. As a consequent of the main result, we obtain a new implicit function theorem in a general infinite-dimensional framework. Furthermore, the achieved results are proved to be instrumental in deriving efficient conditions for (lower and upper) semicontinuity, continuity, calmness, and Lipschitz continuity of optimal value functions associated with parameterized generalized equations in both finite and infinite dimensions. 

The obtained results call for further developments. Among them we mention establishing conditions, which ensure continuity and other regularity properties of solutions to parameterized generalized equations and implicit functions. Another direction of our future research addresses applications of the achieved results to the classes of parametric optimization problems discussed in Section~\ref{intro} that are reduced by parameterized generalized equations with both set-valued terms.

\end{document}